\documentclass[a4paper]{amsart}

\usepackage{amscd}
\usepackage{amsmath}
\usepackage{amssymb}
\usepackage{amsmath,amsthm,amssymb,amsfonts}
\usepackage{mathrsfs}

\author{Sylvie Monniaux}

\address{LATP UMR 6632 - Case cour A - Facult\'e des sciences de Saint-J\'er\^ome - Universit\'e Paul C\'ezanne (Aix-Marseille 3) - 13397 Marseille C\'edex 20 - France}

\email{sylvie.monniaux@univ.u-3mrs.fr}

\title[Navier-Stokes in arbitrary domains]{Navier-Stokes equations in arbitrary domains : \\ the Fujita-Kato scheme}


\def\C{\mathbb C}
\def\R{\mathbb R}
\def\N{\mathbb N}

\def\P{\mathbb P}

\def\H{\mathcal H}
\def\V{\mathcal V}
\def\G{\mathcal G}
\def\D{\mathcal D}
\def\E{\mathcal E}
\def\div{\mbox{\rm div}}

\newcommand{\Con}{\ensuremath{\mathscr{C}}}
\renewcommand{\L}{\ensuremath{\mathscr{L}}}
\newcommand{\Drond}{\ensuremath{\mathscr{D}}}

\theoremstyle{plain}
\newtheorem{proposition}{Proposition}[section]
\newtheorem{theorem}[proposition]{Theorem}

\theoremstyle{definition}
\newtheorem{definition}[proposition]{Definition}

\theoremstyle{remark}
\newtheorem{remark}[proposition]{Remark}

\numberwithin{equation}{section}

\begin{document}

\begin{abstract}
Navier-Stokes equations are investigated in a functional setting in 3D open sets $\Omega$, bounded or not, without assuming any regularity of the boundary $\partial\Omega$. The main idea is to find a correct definition of the Stokes operator in a suitable Hilbert space of divergence-free vectors and apply the Fujita-Kato method, a fixed point procedure, to get a local strong solution.
\end{abstract} 
 
\maketitle

\section{Introduction}

Since the pioneering work by Leray \cite{Le34} in 1934, there have been several studies on solutions of Navier-Stokes equations 
$$
(NS) \ \left\{
\begin{array}{rclcl}
\frac{\partial u}{\partial t} - \Delta u +\nabla\pi + (u\cdot\nabla)u &= &0&\mbox{in}& ]0,T[ \times \Omega,\\
\div \ u &=& 0 &\mbox{in}& ]0,T[ \times \Omega,\\
u&=&0& \mbox{on}& ]0,T[ \times \partial\Omega,\\
u(0) &=& u_0 &\mbox{in}& \Omega.
\end{array}
\right.
$$
Fujita and Kato \cite{FuKa64} in 1964 gave a method to construct so called mild solutions in smooth domains $\Omega$, producing local (in time) smooth solutions of $(NS)$ in a Hilbert space setting. These solutions are global in time if the initial value $u_0$ is small enough in a certain sense. The case of non smooth domains has been studied by Deuring and von Wahl \cite{DeVW95} in 1995 where they considered domains $\Omega\subset \R^3$ with Lipschitz boundary $\partial \Omega$. They found local smooth solutions using results contained in Shen's PhD thesis \cite{Sh91}. Their method does not cover the critical space case as in \cite{FuKa64}. One of the difficulty there was to understand the Stokes operator, and in particular its domain of definition.

In Section 2, we give a ``universal" definition of the Stokes operator, for any domain $\Omega \subset \R^3$ (Defintion \ref{stokes}). In Section 3, we construct a mild solution of $(NS)$ with a method similar to Fujita-Kato's \cite{FuKa64} (Theorem \ref{mildsolutions}) for initial values $u_0$ in the critical space $D(A^{\frac 14})$. We show in Section 4 that this mild solution is a strong solution, $i.e.$ $(NS)$ is satisfied almost everywhere.

\section{The Stokes operator}

Let $\Omega$ be an open set in $\R^3$. The space 
$$
L^2(\Omega)^3 = \{ u = (u_1, u_2, u_3) ; u_i \in L^2(\Omega), \ i=1,2,3\}
$$ 
endowed with the scalar product 
$$
\langle u,v\rangle = \int_\Omega u \cdot \overline{v} = \sum_{i= 1}^3 \int_\Omega u_i \  \overline{v_i}
$$ 
is a Hilbert space. Define
$$
\G = \{ \nabla p ; p \in L_{loc}^2(\Omega) \mbox{ and }\nabla p \in L^2(\Omega)^3\} ;
$$  
the set $\G $ is a closed subspace of $L^2(\Omega)^3$. Let 
$$
\H = \G^\bot = \left\{ u \in L^2(\Omega)^3 ; \langle u , \nabla p \rangle = 0 , \ \forall p \in H^1(\Omega)\right\}.
$$
The space $\H$, endowed with the scalar product $\langle \cdot, \cdot \rangle$ is a Hilbert space. We have the following Hodge decomposition
$$
L^2(\Omega)^3 = \H \stackrel{\bot}{\oplus} \G.
$$
We denote by $\P$ the projection from $L^2(\Omega)^3$ onto $\H$ : $\P$ is the usual Helmoltz projection. We denote by $J$ the canonical injection $\H \hookrightarrow L^2(\Omega)^3$ : $J' = \P$ ($J'$ beeing the adjoint of $J$) and $\P J$ is the identity on $\H$. Let now $\Drond(\Omega)^3 =  \Con_c^{\infty}(\Omega)^3$ and 
$$
\D = \{u \in \Drond(\Omega)^3 ; \div u = 0\}.
$$
It is clear that $\D$ is a closed subspace of $\Drond(\Omega)^3$. We denote by $J_0 : \D \hookrightarrow \Drond(\Omega)^3$ the canonical injection : $J_0 \subset J$. Let $\P_1$ be the adjoint of $J_0$ : $\P_1 = J_0' : \Drond'(\Omega)^3 \to \D'$. We have $\P_1 \subset \P$. The following theorem characterizes the elements in $\ker \P_1$.

\begin{theorem}[de Rahm] \label{deRahm}
Let $T \in \Drond'(\Omega)^3$ such that $\P_1T = 0$ in $\D'$. Then there exists $S \in (\Con_c^\infty(\Omega))'$ such that $T = \nabla S$. Conversely, if $T = \nabla S$ with $S \in (\Con_c^\infty(\Omega))'$, then $\P_1 T = 0$ in $\D'$.
\end{theorem}

We denote by $H^1_0(\Omega)^3$ the closure of $\Drond(\Omega)^3$ with respect to the scalar product $(u,v) \mapsto \langle u, v \rangle_1 = \langle u,v \rangle + \sum_{i=1}^3 \langle \partial_i u , \partial_i v \rangle$. By Sobolev embeddings, we have $H^1_0(\Omega)^3 \hookrightarrow L^6(\Omega)^3$. Define 
$$
\V = \H \cap H^1_0(\Omega)^3.
$$
The space $\V$ is a closed subspace of $H^1_0(\Omega)^3$ ; endowed with the scalar product $\langle \cdot, \cdot \rangle_1$, $\V$ is a Hilbert space. The canonical injection $\tilde{J} : \V \hookrightarrow H^1_0(\Omega)^3$ is the restriction of $J$ to $\V$. Let $H^{-1}(\Omega)^3 = (H^1_0(\Omega)^3)'$ ; $\P_1$ maps $H^{-1}(\Omega)^3$ to $\V'$ : the restriction of $\P_1$ to $H^{-1}(\Omega)^3$ is $\tilde{\P}$, the adjoint of $\tilde{J}$. On $\V \times \V$ we define now the form $a$ by $\displaystyle{a(u,v) = \sum_{i=1}^3 \langle \partial_i \tilde{J}u,  \partial_i \tilde{J}v \rangle}$~: $a$ is a bilinear, symmetric, $\delta + a$ is a coercive form on $\V \times \V$ for all $\delta >0$, then defines a bounded self-adjoint operator $A_0 : \V \to \V'$ by $(A_0u)(v) = a(u,v)$ with $\delta +A_0$ invertible for all $\delta >0$.

\begin{proposition}
For all $u \in \V$, $A_0 u = \tilde{\P} (-\Delta_D^\Omega) \tilde{J}u$, where $\Delta_D^\Omega$ denotes the Dirichlet-Laplacian on $H^1_0(\Omega)^3$.
\end{proposition}

\begin{proof}
For all $u,v \in \V$, we have
\begin{eqnarray*}
(A_0u)(v) &\stackrel{(1)}{=}& a(u,v) 
\ \stackrel{(2)}{=}\ \sum_{i=1}^3 \langle \partial_i \tilde{J}u,  \partial_i \tilde{J}v \rangle \\
&\stackrel{(3)}{=}& \langle (-\Delta_D^\Omega)\tilde{J} u , \tilde{J}v \rangle_{H^{-1},H^1_0}\\
&\stackrel{(4)}{=}& \langle \tilde{\P}(-\Delta_D^\Omega) \tilde{J}u , v \rangle_{\V',\V}.
\end{eqnarray*}
The first two equalities come from the definition of $A_0$ and $a$. The third equality comes from the definition of the Dirichlet-Laplacian on $H^1_0(\Omega)^3$ and the fact that for $v\in\V$, $\tilde{J}v = v$. The last equality is due to $\tilde{J}' \varphi = \tilde{\P} \varphi$ in $\V'$ for all $\varphi \in H^{-1}(\Omega)^3$. This shows that $A_0u$ and $\tilde{\P}(-\Delta_D^\Omega) \tilde{J} u$ are two continuous linear forms on $\V$ which co\"\i ncide on $\V$, they are then equal.
\end{proof}

\begin{definition} \label{stokes}
The operator $A$ defined on its domain $D(A) = \{ u\in \V ; A_0 u \in \H\}$ by $Au = A_0u$ is called the Stokes operator. 
\end{definition}

\begin{theorem}
The Stokes operator is self-adjoint in $\H$, generates an analytic semigroup $(e^{-tA})_{t\ge 0}$, $D(A^{\frac 12}) = \V$ and satisfies
\begin{eqnarray*}
D(A) &= &\{ u \in \V \ ; \ \exists \pi \in (\C_c^\infty(\Omega))' : \nabla\pi \in H^{-1}(\Omega) \mbox{ and } -\Delta u + \nabla \pi \in \H\} \\
Au &=& -\Delta u + \nabla \pi.
\end{eqnarray*}
\end{theorem}

\begin{remark}
Since $H^1_0(\Omega)^3 \hookrightarrow L^6(\Omega)^3$, it is clear by interpolation and dualization that $\P_1$ maps $L^p(\Omega)^3$ to $D(A^s)'$ for $\frac 65\le p\le 2$, $0\le s \le \frac 12$ and $s= -\frac 34 +\frac{3}{2p}$. Since $A$ is self-adjoint, one has $(\delta +A_0)^{-s}D(A^s)' =\{ (\delta +A_0)^{-s}u; u \in D(A^s)'\}= \H$. In particular, $(\delta +A_0)^{-\frac 14}\P_1$ maps $L^{\frac 32}(\Omega)^3$ into $\H$.
\end{remark}

\section{Mild solution to the Navier-Stokes system}

Let $T>0$. 

Define the following Banach space
\begin{eqnarray*}
\E_T &=& \left\{u \in \Con([0,T] ; D(A^{\frac 14})\cap \Con^1(]0,T] ; D(A^{\frac 14})) \right.\\
&&\left. \mbox{such that } \ \sup_{0<s < T} \|s^{\frac 14} A^{\frac 12}u(s) \|_\H + \sup_{0<s <T} \|s A^{\frac 14}u'(s) \|_\H < \infty \right\}
\end{eqnarray*}
endowed with the norm
$$
\|u\|_{\E_T } = \sup_{0<s<T} \|A^{\frac 14} u(s) \|_\H + \sup_{0<s<T} \|s^{\frac 14} A^{\frac 12}u(s) \|_\H + \sup_{0<s<T} \|s A^{\frac 14}u'(s) \|_\H .
$$
Let $\alpha$ be defined by $\alpha(t) = e^{-tA} u_0$ where $u_0 \in D(A^{\frac 14})$. Then $\alpha \in \E_T$. Indeed, it is clear that $\alpha \in \Con([0,T] ; D(A^{\frac 14}))$. We also have that $t^{\frac 14} A^{\frac 12} \alpha(t) = t^{\frac 14} A^{\frac 14} e^{-tA} A^{\frac 14} u_0$ is bounded on $(0,T)$ since $(e^{-tA})_{t \ge 0}$ is an analytic semigroup. Moreover, one has $\alpha'(t) = -A e^{-tA} u_0$ which yields to $t A^{\frac 14} \alpha'(t) = - tAe^{-tA} A^{\frac 14}u_0$ continuous on $]0,T]$, bounded in $\H$. For $u,v \in \E_T$, we define now
$$
\Phi (u,v)( t ) = \int_0^t e^{-(t-s) A} \textstyle{(-\frac 12 \P_1)} ((u(s) \cdot \nabla) v(s) + (v(s) \cdot \nabla) u(s)) ds , \quad 0<t<T.
$$

\begin{proposition}
The transform $\Phi$ is bilinear, symmetric, continuous from ${\mathcal E}_T \times {\mathcal E}_T$ to ${\mathcal E}_T$ and the norm of $\Phi$ is independent of $T$.
\end{proposition}

\begin{proof}
The fact that $\Phi$ is bilinear and symmetric is clear. Moreover, $\Phi(u,v) = e^{-\cdot A} * f$, where $f$ is defined by 
$$
f(s) =  \textstyle{(-\frac 12 \P_1)}((u(s) \cdot \nabla) v(s) + (v(s) \cdot \nabla) u(s)), \quad s \in [0,T].
$$
For $u,v \in {\mathcal E}_T$, it is clear that $(u(s) \cdot \nabla) v(s) + (v(s) \cdot \nabla) u(s) \in L^{\frac 32}(\Omega)^3$ and therefore $(\delta +A_0)^{- \frac 14} f(s) \in \H$ with $\displaystyle{\sup_{0<s<T} s^{\frac 12}\|(\delta +A_0)^{- \frac 14}  f(s)\|_\H \le c \|u\|_{\E_T} \|v\|_{\E_T}}$. We have then
$$
\Phi(u,v) = e^{-\cdot A} * f = (\delta +A)^{\frac 14}e^{-\cdot A} * ((\delta +A_0)^{- \frac 14} f)
$$ 
and therefore
\begin{eqnarray*}
\|A^{\frac 14} \Phi(u,v)(t) \|_\H 
&\le& \int_0^t \|A^{\frac 14} (\delta +A)^{\frac 14} e^{-(t-s)A}\|_{\L(\H)}\|(\delta +A_0)^{- \frac 14}f(s)\|_\H ds \\
&\le& c \left(\int_0^t \frac{1}{\sqrt{t-s}} \frac{1}{\sqrt{s}} \ ds \right) \|u\|_{\E_T}  \|v\|_{\E_T} \\
&\le& c \left(\int_0^1 \frac{1}{\sqrt{1-\sigma}} \frac{1}{\sqrt{\sigma}} \ d\sigma \right) \|u\|_{\E_T}  \|v\|_{\E_T} \\
&\le& c \|u\|_{\E_T}  \|v\|_{\E_T}.
\end{eqnarray*}
Continuity with respect to $t \in [0,T]$ of $t \mapsto A^{\frac 14} \Phi(u,v)(t)$ is clear once we have proved the boundedness. We also have
\begin{eqnarray*}
\|A^{\frac 12} \Phi(u,v)(t) \|_\H
&\le&  \int_0^t \|A^{\frac 12} (\delta +A)^{\frac 14}e^{-(t-s)A}\|_{\L(\H)}\|(\delta +A_0)^{- \frac 14}f(s)\|_\H ds \\
&\le& c \left(\int_0^t \frac{1}{(t-s)^{\frac 34}} \frac{1}{\sqrt{s}} \ ds \right) \|u\|_{\E_T}  \|v\|_{\E_T} \\
&\le& c t^{-\frac 14} \left(\int_0^1 \frac{1}{(1-\sigma)^{\frac 34}} \frac{1}{\sqrt{\sigma}} \ d\sigma \right) \|u\|_{\E_T}  \|v\|_{\E_T} \\
&\le& c t^{-\frac 14} \|u\|_{\E_T}  \|v\|_{\E_T}.
\end{eqnarray*}
Continuity with respect to $t \in ]0,T]$ is clear once we have proved the boundedness. To prove the last part of the norm of $\Phi(u,v)$ in $\E_T$, we have for $s \in ]0,T[$
$$
f'(s) =  \textstyle{(-\frac 12 \P_1)}((u'(s) \cdot \nabla) v(s) + (u(s) \cdot \nabla) v'(s) + (v'(s) \cdot \nabla) u(s) + (v(s) \cdot \nabla) u'(s))
$$
and therefore
$$
\sup_{0<s<T} \|s^{\frac 54} (\delta +A_0)^{-\frac 12} f'(s)\|_\H \le c \|u\|_{\E_T}  \|v\|_{\E_T}.
$$
We have
$$
\Phi (u,v) (t) = \int_0^{\frac t2} e^{-sA} f(t-s) ds + \int_0^{\frac t2} e^{-(t-s)A} f(s) ds \quad t \in ]0,T[,
$$
and therefore
\begin{eqnarray*}
\Phi(u,v)'(t) &=& e^{-\frac t2 A} \textstyle{f (\frac t2)} + 
\displaystyle{\int_0^{\frac t2} (\delta +A)^{\frac 12} e^{-sA} (\delta +A_0)^{-\frac 12} f'(t-s) ds} \\
&& + \int_0^{\frac t2} -A (\delta +A)^{\frac 14} e^{-(t-s)A} (\delta +A_0)^{-\frac 14}f(s) ds,
\end{eqnarray*}
which yields
\begin{eqnarray*}
\|A^{\frac 14} \Phi(u,v)'(t)\|_\H &\le&  \frac{c}{\sqrt{t}}\textstyle{ \left\|(\delta +A_0)^{-\frac 14} f(\frac t2)\right\|_\H} + c \displaystyle{\left(\int_0^{\frac t2} \frac{1}{s^{\frac 12}} \frac{1}{(t-s)^{\frac 54}} ds\right) \|u\|_{\E_T}  \|v\|_{\E_T} }\\
&& + c \left(\int_0^{\frac t2} \frac{1}{(t-s)^{\frac 54}} \frac{1}{s^{\frac 12}} ds\right) \|u\|_{\E_T}  \|v\|_{\E_T}\\
&\le& \frac{c}{t} \left(\int_0^{\frac 12} \frac{d\sigma}{(1-\sigma)^{\frac 54} \sigma^{\frac 12}} \right) \|u\|_{\E_T}  \|v\|_{\E_T}.
\end{eqnarray*}
This last inequality ensures that $\Phi(u,v) \in \E_T$ whenever $u,v \in \E_T$.
\end{proof}

\begin{theorem} \label{mildsolutions}
For all $u_0 \in D(A^{\frac 14})$, there exists $T>0$ such that there exists a unique $u \in {\mathcal E}_T$ solution of $u = \alpha + \Phi(u,u)$ on $[0,T]$. This function $u$ is called the mild solution to the Navier-Stokes system.
\end{theorem}

\begin{proof}
Let $T>0$. Since $\Phi : \E_T \times \E_T \to \E_T$ is bilinear continuous, it suffices to apply Picard fixed point theorem, as in \cite{FuKa64}. The sequence in $\E_T$ $(v_n)_{n\in\N}$ defined by $v_0=\alpha$ as first term and 
$$
v_{n+1} = \alpha + \Phi(v_n,v_n), \quad n\in \N
$$
converges to the unique solution $u \in \E_T$ of $u = \alpha + \Phi(u,u)$ provided $\|A^{\frac 14}u_0\|_\H$ is small enough ($\|\alpha\|_{\E_T} < \frac{1}{4 \|\Phi\|_{\L(\E_T \times \E_T ; \E_T)}}$). In the case where $\|A^{\frac 14}u_0\|_\H$ is not small (that is, if $\|\alpha\|_{\E_T} \ge \frac{1}{4 \|\Phi\|_{\L(\E_T \times \E_T ; \E_T)}}$) then for $\varepsilon >0$, there exists $u_{0,\varepsilon} \in D(A)$ such that $\|A^{\frac 14} (u_0-u_{0,\varepsilon})\|_\H \le \varepsilon$. If we take as initial value $u_{0, \varepsilon} \in D(A)$, we have
$$
\|\alpha_{\varepsilon}\|_{\E_T} \le c T^{\frac 34} \|Au_{0,\varepsilon}\|_\H \xrightarrow[T\to 0]{} 0.
$$
Therefore, we can find $T>0$ such that $\|\alpha\|_{\E_T} < \frac{1}{4 \|\Phi\|_{\L(\E_T \times \E_T ; \E_T)}}$.
\end{proof}

\section{Strong solutions}

Let $u$ be the mild solution to the Navier-Stokes system. We show in this section that $u$ in fact satisfies the equations of the Navier-Stokes system in an $L^p-$sense (for a suitable $p$). To begin with, we know that $u \in {\mathcal E}_T$ and satisfies 
$$
u = \alpha + \Phi(u,u) = \alpha + e^{-\cdot A} * \varphi(u),
$$
where $\varphi(u) =  - \P_1((u \cdot \nabla) u)$ and we have $\|t^{\frac 12}(u(t)\cdot\nabla)u(t)\|_{\frac 32} \le c \|u\|_{\E_T}^2$. Therefore, we get 
\begin{equation} \label{inival}
u(0) = \alpha(0) = u_0, 
\end{equation}
\begin{equation} \label{div}
\div u(t) = 0 \mbox{ in the } L^2-\mbox{sense for } t\in]0,T[,
\end{equation}
and 
$$
u' + Au = f \quad \mbox{ in } \Con(]0,T[ ; \V'),
$$
which means that for all $t\in]0,T[$, 
$$
\P_1(u'(t) -\Delta_D^\Omega u(t) + (u(t)\cdot\nabla)u(t) ) = 0.
$$ 
Then, by Theorem~\ref{deRahm}, there exists $(-\pi)(t) \in (\Con_c^\infty(\Omega))'$ such that $\nabla \pi(t) \in H^{-1}(\Omega)^3$ and
\begin{equation}\label{equa}
\nabla (-\pi)(t) = u'(t) -\Delta_D^\Omega u(t) + (u(t)\cdot\nabla)u(t) 
\end{equation}
and we have for $0<t<T$
$$
-\Delta_D^\Omega u(t) +\nabla \pi(t) = -u'(t) - (u(t)\cdot\nabla)u(t)  \in L^3(\Omega)^3 +  L^{\frac 32}(\Omega)^3 .
$$
The equation \eqref{equa}, together with \eqref{inival} and \eqref{div}, give the usual Navier-Stokes equations which are fulfilled in a strong sense ($a.e.$) where we consider the expression $-\Delta u +\nabla \pi$ undecoupled.

\end{document}